\documentclass{article}
\usepackage{amsmath,amssymb,epsf}
\usepackage{theorem}

\newtheorem{theorem}{Theorem}

\newtheorem{dfn}{Definition}
\newtheorem{cor}{Corollary}

\newenvironment{example}{\noindent\textbf{Examples}}{}
\newenvironment{prf}{\noindent\textsc{Proof}}{\hfill\ensuremath{\blacksquare}}
\newenvironment{remark}{\noindent\textbf{Remark}}{}
\newenvironment{ack}{\noindent\textbf{Acknowledgments}}{}
\begin{document}
\title{A Multiplication Rule for\\ the Descent Algebra
of Type
$D$} 
\author{N. Bergeron and S.J. van Willigenburg\\Department of Mathematics and 
Statistics,\\York University, 4700 Keele St,\\North York, 
ON, M3J 1P3, CANADA.} 
\maketitle 
\begin{abstract} Here we give an interpretation of Solomon's 
rule  for multiplication in the descent algebra of Coxeter
groups of type $D$, $\Sigma D_n$. We describe an ideal
${\cal I}$ such that 
$\Sigma D_n /{\cal I}$ is isomorphic to the descent
algebra of the hyperoctahedral group,
$\Sigma B_{n-2}$.\end{abstract}

\section{Introduction}

Given a Coxeter group, $W$, we can construct an algebra - \textit{the descent
algebra} - which is a sub-algebra of the group algebra $\mathbb{Q}[W]$. These
were introduced in 1976 by Louis Solomon 
\cite{solomon-mackey}. A revival of interest in this
area began in the 80's when applications were found
for an interpretation of the rule  for multiplying
together basis elements of the descent algebra of the
symmetric group, for example
\cite{garsia-reutenauer}, \cite{garsia-remmel}. Since this 
interpretation involved matrices,
 we shall call it the ``matrix interpretation'' from here on. This  matrix interpretation provided
the key to many advances in the subject (for instance  \cite{atkinson-solomon},
\cite{bergeron-gr}, \cite{bergberg-ht}) including an analogous matrix interpretation by
Fran\c{c}ois and Nantel Bergeron  for the descent algebra of the hyperoctahedral group, 
\cite{bergeron-bergeron}.

Until now, there has been little success in developing such an interpretation for the Coxeter groups
of type $D$. However, in this paper we shall give the matrix interpretation for this remaining
Coxeter family, after defining the Coxeter groups of type $D$, and their corresponding descent
algebra.

The $n$-th Coxeter group of type $D$, $D_n$, is the group
acting on the set $$\{-n, \ldots ,-1,1,\ldots ,n\}$$ whose
Coxeter generators are the set
$S=\{s_{1'},s_1,s_2,\ldots ,s_{n-1}\}$, where 
$s_i$ is the product of transpositions $(-i\!
-\! 1,\, -i)(i,i+1)$ for
$i=1,2,\ldots ,n-1$, and
$s_{1'}$ is the product of transpositions $(-2,1)(-1,2)$.
The  relations are given by the following diagram:

\begin{figure}[htbp]
\centerline{\epsffile{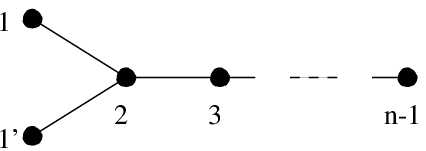}}
\end{figure}

where an edge between distinct nodes $i$ and $j$ gives us the relation $(s_is_j)^3=1$, and
no edge gives $(s_is_j)^2=1$, and $(s_i)^2=1$.

Solomon proved that if $J$ is a subset of  $S$,  
$W_J$ is the subgroup generated by $J$,  
$X_J$ ($X_J^{-1}$) is the unique set of minimal
length left (right) coset representatives of $W_J$,
and ${\cal X}_J$ is the formal sum of the elements
in $X_J$ then for $J,K,L\subseteq S$

\[{\cal X}_J{\cal X}_K=\sum _L a_{JKL}{\cal X}_L\]

where $a_{JKL}$ is the number of elements $x\in X_J^{-1}\cap X_K$ such that
$x^{-1}Jx\cap K=L$. Hence, the set of all ${\cal X}_J$'s form a basis for an
algebra - the descent algebra of  $D_n$, $\Sigma D_n$. Our interpretation of
this multiplication rule uses this basis, but for ease of computation, we use a
different notation.

We define a composition, $q$, of an integer, $n$, to be an ordered list  
$[q_1,q_2,\ldots ,q_k]$ of positive integers whose sum is $n$, and shall write
$q\vDash n$ to denote this. We shall call the integers $q_1,q_2,\ldots ,q_k$ the 
\textit{components} of
$q$.

 There exists a natural bijection between the subsets of $S$ and the disjoint union, 
$C(n)$, of the sets \( C_{<n}=\{q|q\vDash m,m\leq
	n-2\}\), \( C_1=\{q|q\vDash n,q_1=1 \}\) 
\( C_n=\{q|q\vDash n,q_1\geq 2 \}\) and 
  \( C_n'=\{q|q\vDash n,q_1\geq
2\}\). Note that   $C_n$ and
$C_n'$ are two  copies of the same set. Let $q\in C(n)$ such that
$q\vDash m\leq n$, then  	the 	subset 	 
corresponding to $q$ is
 \begin{enumerate}
\item \(\{s_{q_0}, s_{q_0+q_1},\ldots , s_{q_0+\ldots
+q_{(k-1)}}\}\) if \(q\in\)$C_{<n}$,
\item \(\{s_{1'}, s_{1}, s_{1+q_2},\ldots , s_{1+q_2+\ldots +q_{(k-1)}}\}\)  if \(q\in\)\( C_1\),
\item \(\{s_{1'},s_{ q_1},\ldots ,  s_{q_1+\ldots
+q_{(k-1)}}\}\)  if \(q\in\)\( C_n\), 
\item \(\{s_{1}, s_{q_1},\ldots , s_{q_1+\ldots +q_{(k-1)}}\}\) if \(q\in\)\( C_n'\),
 
\end{enumerate}
where \(q_0=n-m.\)

\begin{remark}
The   step of corresponding a set, $J$, containing
$s_{1'}$ ($s_1$) with a composition, $q$, in $C_n$
($C_n'$) is because we shall later relate $q$ to the
complement of $J$.
\end{remark}

\section{The Matrix Interpretation, and  Results}

If $J^c$ is the complement of $J$ in $S$, then we let
$B_q={\cal X}_{ J^c}$ where $q$ is the composition in $C(n)$ that corresponds to
$J$ by the
 above bijection. The matrix interpretation of Solomon's multiplication rule can now be described as follows.

Consider the template with the following form
\[\begin{pmatrix}
z_{00}&z_{01}&z_{02}&\ldots &z_{0v}\\
&y_{11}&y_{12}&\ldots &y_{1v}\\
z_{10}&z_{11}&z_{12}&\ldots &z_{1v}\\
\vdots &\vdots &\vdots &\ddots &\vdots\\
&y_{u1}&y_{u2}&\ldots &y_{uv}\\
z_{u0}&z_{u1}&z_{u2}&\ldots &z_{uv}\\
\end{pmatrix}\]
Note that the $y$-lines do not have entries in column 0.
We say a template is a ``filled template'' 
if all entries in a  template are non-negative integers.

\begin{dfn} Let $\boldsymbol{t}$ be a filled template.   We
define the \emph{border-sum},  ${\cal B}(\boldsymbol{t})$,
of $\boldsymbol{t}$ to be the sum
\[z_{00}+\sum _{i=1}^{u}z_{i0}+\sum _{j=1}^{v}z_{0j}\]
and the\emph{ y-sum}, ${\cal Y}(\boldsymbol{t})$, to be \(\sum_{i, j}
y_{ij}\). The reading word of $\boldsymbol{t}$, $r(\boldsymbol{t})$,  is given
by {\small \[[z_{01},z_{02},\ldots ,z_{0v},y_{1v},\ldots
,y_{12},y_{11},z_{10},z_{11},z_{12},\ldots  ,z_{1v},\ldots ,z_{u0},z_{u1},z_{u2},
\ldots ,z_{uv}]\]} with zero entries omitted, unless $z_{00}=1$, in which case $r(\boldsymbol{t})$
is given by
 {\small \[[1,z_{01},z_{02},\ldots ,z_{0v},y_{1v},\ldots
,y_{12},y_{11},z_{10},z_{11},z_{12},\ldots  ,z_{1v},\ldots ,z_{u0},z_{u1},z_{u2},
\ldots ,z_{uv}]\]}with zero entries
omitted.\end{dfn}

If \(p\), and
\(q\) are compositions in $C (n)$  such that $p\vDash l\leq n$, and $q\vDash m\leq n$,  then let 
\(Z(p,q)\) be the set of  filled templates, $\boldsymbol{t}$, such that
\begin{enumerate}
\item \(z_{0j}+\sum _{i\neq 0}(y_{ij}+z_{ij})=p _{j},\ j\neq 0\),
\item $\sum _i z_{i0}=n-l$,
\item \(z_{i0}+\sum _{j\neq 0}(y_{ij}+z_{ij})=q _{i},\ i\neq 0\),
\item $\sum _{j} z_{0j}=n-m$,
\item If ${\cal B}(\boldsymbol{t})$=0, ${\cal Y}(\boldsymbol{t})$ is odd if
 \begin{enumerate}\item $p\in C _1\cup C _n$
and $q\in C _n'$, or 
\item $p\in C _n'$ and $q\in C _1\cup C _n$.\end{enumerate}
 Otherwise ${\cal Y}(\boldsymbol{t})$ is even.
\end{enumerate}

We are now ready to state our matrix interpretation. To distinguish between those compositions
belonging to $C_n$ and those belonging to $C_n'$, we shall write $q'$ when $q\vDash n$ and $q\in
C_n'$.

\begin{theorem}
Let $p,q\in C (n)$. For any filled template
$\boldsymbol{t}$, let
$r(\boldsymbol{t})=[r_1(\boldsymbol{t}), r_2(\boldsymbol{t}), 
\ldots ]$.
  Then,
\[B_p B_q=\sum _{\boldsymbol{t}\in Z(p,q)}
\tilde{B}_{r(\boldsymbol{t})}\] where
$\tilde{B}_{r(\boldsymbol{t})}$ satisfies the following.

\begin{enumerate}
\item If $q\in C _1$, then  
$\tilde{B}_{r(\boldsymbol{t})}=B_{r(\boldsymbol{t})}
$.\label{ono-p}
\item If $q\in C _n$, then  
$\tilde{B}_{r(\boldsymbol{t})}=B_{r(\boldsymbol{t})}
$.\label{no-p}
\item  If $q\in C _n'$, then if 
$r_1(\boldsymbol{t})=1$ then $\tilde{B}_{r(\boldsymbol{t})}=B_{r(\boldsymbol{t})} $,
otherwise
$\tilde{B}_{r(\boldsymbol{t})}=B_{r(\boldsymbol{t})'}$.
\label{rt-up}
\item  If $q\in C _{<n}$, then \label{ty-pe1}
\begin{enumerate}
\item If  $r_1(\boldsymbol{t})\geq 2$, $p\in C _1\cup C _n$ and ${\cal
Y}(\boldsymbol{t})$ is odd, or $r_1(\boldsymbol{t})\geq 2$,
  $p\in C _n'$ and ${\cal Y}(\boldsymbol{t})$ is even,
then 
$\tilde{B}_{r(\boldsymbol{t})}=B_{r(\boldsymbol{t})'}$.
\label{on-ep}
\item If $p\in C _{<n}$ and $z_{00}=0$, then if
$r_1(\boldsymbol{t}) =1$, then $\tilde{B}_{r(\boldsymbol{t})}=2B_{r(\boldsymbol{t})}
$, otherwise
$\tilde{B}_{r(\boldsymbol{t})}=B_{r(\boldsymbol{t})}
+B_{r(\boldsymbol{t}) '}$.  \label{two-u}
\item Otherwise
$\tilde{B}_{r(\boldsymbol{t})}=B_{r(\boldsymbol{t})}$.
\end{enumerate} 
\end{enumerate}
\label{d-mult}
\end{theorem}

A rigorous proof of this theorem can be obtained through a variety of methods. One is to use
shuffle products in a way similar to that seen in \cite{garsia-remmel}, or sketched in
\cite{bergeron-bergeron}, to prove the analogous interpretations for the descent algebras of the
Coxeter groups of types $A$ and $B$, respectively. Alternatively, this theorem can be proved using
the general framework suggested in \cite{vanwilli-justgar}. Indeed this framework inspired Theorem
~\ref{d-mult}, and a proof in this vein can be found in \cite{vanwilli-thesis}.

Here, however, we wish to emphasize that it is the formulation of Theorem ~\ref{d-mult} that is the
most difficult stage. Once this has been achieved, a proof can be derived by the diligent reader,
with or without the use of the above references, or found in \cite{vanwilli-thesis}.
Therefore we feel it would be more beneficial to replace the proof with a collection of
illuminating examples.

\begin{example}
To illustrate our rule we shall work in $\Sigma D_4$.  Each
example, $B_pB_q$,  shall  consist of
$Z(p,q)$, and the resulting summands it generates according
to the rule.
\begin{enumerate}
\item $B_{[4]}B_{[1,3]}$.
\[\begin{pmatrix}
0&0\\
&1\\
0&0\\
&3\\
0&0
\end{pmatrix}
\begin{pmatrix}
0&0\\
&0\\
0&1\\
&0\\
0&3
\end{pmatrix}
\begin{pmatrix}
0&0\\
&0\\
0&1\\
&2\\
0&1
\end{pmatrix}
\begin{pmatrix}
0&0\\
&1\\
0&0\\
&1\\
0&2
\end{pmatrix}\]
\[B_{[4]}B_{[1,3]}=2B_{[1,3]}+B_{[1,2,1]}+B_{[1,1,2]}\]\label{b31-b4}
\item  $B_{[3,1]'}B_{[4]}$.
\[\begin{pmatrix}
0&0&0\\
&3&0\\
0&0&1
\end{pmatrix}
\begin{pmatrix}
0&0&0\\
&0&1\\
0&3&0
\end{pmatrix}
\begin{pmatrix}
0&0&0\\
&2&1\\
0&1&0
\end{pmatrix}
\begin{pmatrix}
0&0&0\\
&1&0\\
0&2&1
\end{pmatrix}\]
\[B_{[3,1]'}B_{[4]}=B_{[3,1]}+B_{[1,3]}+2B_{[1,2,1]}\]\label{b31p-b4}
\item  $B_{[2,2]'}B_{[4]'}$.
\[\begin{pmatrix}
0&0&0\\
&2&2\\
0&0&0
\end{pmatrix}
\begin{pmatrix}
0&0&0\\
&2&0\\
0&0&2
\end{pmatrix}
\begin{pmatrix}
0&0&0\\
&0&2\\
0&2&0
\end{pmatrix}
\begin{pmatrix}
0&0&0\\
&0&0\\
0&2&2
\end{pmatrix}
\begin{pmatrix}
0&0&0\\
&1&1\\
0&1&1
\end{pmatrix}
\]
\[B_{[2,2]'}B_{[4]'}=4B_{[2,2]'}+B_{[1,3]}+B_{[1,1,1,1]}\]\label{b22p-b4p}
\item $B_{[4]}B_{[2]}$.
\[\begin{pmatrix}
0&2\\
&2\\
0&0
\end{pmatrix}\quad
\begin{pmatrix}
0&2\\
&0\\
0&2
\end{pmatrix}\quad
\begin{pmatrix}
0&2\\
&1\\
0&1
\end{pmatrix}\]
\[B_{[4]}B_{[2]}=2B_{[2,2]}+B_{[2,1,1]'}\]\label{b4-b2}
\item $B_{[2]}B_{[2]}$.
\[
\begin{pmatrix}
2&0\\
&2\\
0&0
\end{pmatrix}
\begin{pmatrix}
2&0\\
&0\\
0&2
\end{pmatrix}
\begin{pmatrix}
2&0\\
&1\\
0&1
\end{pmatrix}
\begin{pmatrix}
0&2\\
&0\\
2&0
\end{pmatrix}
\begin{pmatrix}
1&1\\
&1\\
1&0
\end{pmatrix}
\begin{pmatrix}
1&1\\
&0\\
1&1
\end{pmatrix}\]
\[B_{[2]}B_{[2]}=2B_{[2]}+B_{[1,1]}+B_{[2,2]}+B_{[2,2]'}+2B_{[1,1,1,1]}\]\label{b2-b2}
\item $B_{[1,1]}B_{[2]}$
\begin{gather*}\begin{pmatrix}
2&0&0\\
&1&1\\
0&0&0
\end{pmatrix}
\begin{pmatrix}
2&0&0\\
&1&0\\
0&0&1
\end{pmatrix}
\begin{pmatrix}
2&0&0\\
&0&1\\
0&1&0
\end{pmatrix}
\begin{pmatrix}
2&0&0\\
&0&0\\
0&1&1
\end{pmatrix}
\begin{pmatrix}
0&1&1\\
&0&0\\
2&0&0
\end{pmatrix}\\
\begin{pmatrix}
1&1&0\\
&0&1\\
1&0&0
\end{pmatrix}
\begin{pmatrix}
1&1&0\\
&0&0\\
1&0&1
\end{pmatrix}
\begin{pmatrix}
1&0&1\\
&1&0\\
1&0&0
\end{pmatrix}
\begin{pmatrix}
1&0&1\\
&0&0\\
1&1&0
\end{pmatrix}\end{gather*}
\[B_{[1,1]}B_{[2]}=4B_{[1,1]}+2B_{[1,1,2]}+4B_{[1,1,1,1]}\]\label{b11-b2}
\end{enumerate}

\begin{remark}
Note, in particular, that these examples illustrate the various conditions given in  Theorem
~\ref{d-mult}.  Examples ~\ref{b31-b4} and ~\ref{b31p-b4} illustrate conditions ~\ref{ono-p} and
~\ref{no-p} respectively, and   the  influence of ${\cal
B}(\boldsymbol{t})=0$ on possible filled templates belonging
to $Z(p,q)$. Example ~\ref{b22p-b4p}  illustrates condition
~\ref{rt-up}, and examples  ~\ref{b4-b2}, ~\ref{b2-b2} and
~\ref{b11-b2} illustrate condition ~\ref{ty-pe1}. More
specifically, examples   ~\ref{b4-b2}, ~\ref{b2-b2} and
~\ref{b11-b2} illustrate respectively conditions
~\ref{on-ep}, ~\ref{two-u} when
$r _1(\boldsymbol{t})\geq 2$, and    ~\ref{two-u} when $r _1(\boldsymbol{t})=1$.
\end{remark}
\end{example}

\begin{cor}
${\cal I}=<B_q|q\in
C_1\cup C_n\cup C_n'>$ is an ideal.
\end{cor}
\begin{prf}
Let $B_p$ be a basis element of $\Sigma D_n$, and $B_q\in {\cal I}$. From our matrix
interpretation it follows that any filled template, $T$, in $Z(p,q)$ or $Z(q,p)$ will
be such that $z_{00}=0$. Therefore $r(T)\vDash n$, that is $B_pB_q,B_qB_p\in {\cal
I}$. The corollary follows immediately by linearity.
\end{prf}

Moreover, we have the following.

\begin{theorem}
Let $B_n$ be the Coxeter group of type $B$, whose  Dynkin
diagram is on $n$ vertices, and let $\Sigma B_n$ be its
associated descent algebra. Then
\[\Sigma B_{n-2}\cong \Sigma D_n/{\cal I}\]
\end{theorem}
\begin{prf}
For clarity, for $q\in$$C_{<n}$, let $B_q^D$ be a basis element of $\Sigma D_n$, and let $B_q^B$ be a basis element of $\Sigma B_{n-2}$.

Note that the set $\{B_q^D\}_{q\in C_{<n}}$ is
a basis for $ \Sigma D_n/{\cal I}$. Hence, let
$p\vDash m_1$, $q\vDash m_2$, $m_1,m_2\leq
n-2$.

By Theorem ~\ref{d-mult}, it follows that in   $\Sigma
D_n/{\cal I}$, the only non-zero term in the product
$B_p^DB_q^D$ are those corresponding to filled templates in
$Z(p,q)$ with $z_{00}\geq 2$. We denote this set of filled
templates by $I(p,q)$. Note that if we subtract 2 from the
$z_{00}$ of any filled template, $T\in I(p,q)$, the reading
word, row sum, and column sum of $T$ are unaffected.
Moreover, if this is performed on all 
$T\in I(p,q)$ the resulting filled templates are precisely
those that arise if we calculate the product $B_p^BB_q^B$ in
$\Sigma B_{n-2}$ (\cite{bergeron-bergeron}). Since this
argument is reversible, the result follows.
\end{prf}

\begin{ack}
The authors are indebted to Michael Atkinson for the opportunity to work together, 
and to him, G\"{o}tz Pfeiffer and Meinolf Geck for many useful discussions.
\end{ack}

\end{document}